\documentstyle[amssymb]{amsart}

\title{On the logical strength of Nash-Williams' theorem on 
transfinite sequences}
\author{Alberto Marcone}
\address{Dipartimento di Matematica, Universit\`a di Torino,
via Carlo Alberto 10, 10123 Torino, Italy}
\email{marcone@@dm.unito.it}
\subjclass{Primary: 03F35. Secondary: 06A07}
\keywords{Better-quasi-ordering, reverse mathematics, minimal bad array}
\thanks{This paper has been completed during the summer of 1994, 
when the author was visiting University of Illinois at
Urbana-Champaign with funding by MURST of Italy and by a 
bilateral CNR-NSF grant.}

\newtheorem{theorem}{Theorem}[section]
\newtheorem{lemma}[theorem]{Lemma}
\newtheorem{corollary}[theorem]{Corollary}
\newtheorem{conjecture}[theorem]{Conjecture}

\theoremstyle{definition}
\newtheorem{definition}[theorem]{Definition}
\newtheorem{remark}[theorem]{Remark}

\renewcommand{\P}[2]{$\Pi_{#2}^{#1}$}
\renewcommand{\S}[2]{$\Sigma_{#2}^{#1}$}
\newcommand{\D}[2]{$\Delta_{#2}^{#1}$}
\newcommand{\bP}[2]{${\boldsymbol\Pi}_{#2}^{#1}$}
\newcommand{\bS}[2]{${\boldsymbol\Sigma}_{#2}^{#1}$}
\newcommand{\NWT}{{\normalshape\bf NWT}}
\newcommand{\GHT}{{\normalshape\bf GHT}}
\newcommand{\Qt}{\tilde Q}
\newcommand{\Bt}{B^*}
\newcommand{\Qom}{Q^{<\omega}}
\newcommand{\RCA}{{\normalshape\bf RCA${}_0$}}
\newcommand{\ACA}{{\normalshape\bf ACA${}_0$}}
\newcommand{\SAC}{{\normalshape\mbox{\bS11-{\bf AC${}_0$}}}}
\newcommand{\ATR}{{\normalshape\bf ATR${}_0$}}
\newcommand{\PCA}{{\normalshape\mbox{\bP11-{\bf CA${}_0$}}}}
\newcommand{\PPCA}{{\normalshape\mbox{\bP12-{\bf CA${}_0$}}}}        
\newcommand{\M}{\frak M}
\newcommand{\Sq}[1]{\left[{#1}\right]^{\omega}}
\newcommand{\sq}[1]{\left[{#1}\right]^{<\omega}}
\newcommand{\set}[2]{\left\{\,{#1}\mid{#2}\,\right\}}
\newcommand{\imply}{\;\longrightarrow\;}
\newcommand{\sse}{\longleftrightarrow}
\newcommand{\lh}{\operatorname{lh}}
\newcommand{\dom}{\operatorname{dom}} 

\newcommand{\base}{\operatorname{base}}

\newcommand{\conc}{{}^\smallfrown}
\newcommand{\sconc}{\conc\!\!}
\newcommand{\N}{\Bbb N}
\newcommand{\init}{\sqsubset}
\newcommand{\initeq}{\sqsubseteq}
\newcommand{\KB}{\leq_{KB}}
\newcommand{\<}{\left\langle}
\renewcommand{\>}{\right\rangle}
\newcommand{\tri}{\vartriangleleft}
\newcommand{\ntri}{\ntriangleleft}
\newcommand{\s}{\sigma}
\newcommand{\pz}{\pi_0}
\newcommand{\po}{\pi_1}
\newcommand{\fom}{f_\omega}
\renewcommand{\H}{\cal H}
\renewcommand{\ss}[2]{\sigma\restriction\Int{#1}{#2}}
\newcommand{\Int}[2]{\left[{#1},{#2}\right)}

\begin{document}
\pagestyle{plain}
\maketitle

\begin{abstract}
We show that Nash-Williams' theorem asserting that the countable 
transfinite sequences of elements of a better-quasi-ordering 
ordered by embeddability form a better-quasi-ordering is provable 
in the subsystem of second order arithmetic \PCA\ but is not 
equivalent to \PCA. We obtain some partial results towards the 
proof of this theorem in the weaker subsystem \ATR\ and we show 
that the minimality lemmas typical of wqo and bqo theory imply 
\PCA\ and hence cannot be used in such a proof. 
\end{abstract} 

The most natural generalization of the notion of well-ordering to 
partial orderings is the concept of well-quasi-ordering or wqo. A 
binary relation $\preceq$ on a set $Q$ is a quasi-ordering if it 
is reflexive and transitive (it is a partial ordering if it is 
also anti-symmetric) and a quasi-ordering is wqo if it contains 
no infinite descending sequences and no infinite anti-chains 
(sets of mutually incomparable elements). Ramsey's theorem 
implies that a quasi-ordering $(Q, \preceq)$ is wqo if and only if for any 
function $f:\N \to Q$ there exist $n$ and $m$ such that $n<m$ and 
$f(n) \preceq f(m)$. For a history of the notion of wqo see 
\cite{Kruskal}. 

Unfortunately the notion of wqo does not enjoy nice closure 
properties: when we apply an operation in some sense infinitary 
to a wqo we cannot expect the resulting quasi-ordering to be wqo. 
This shortcoming prompted Nash-Williams, in a sequence of papers 
culminating in \cite{NW2}, to introduce the notion of
better-quasi-ordering or bqo, i.e.\ wqo of a particularly 
well-behaved kind. Since then a vast theory of bqos has been 
developed and has provided the tools for proving that many 
natural quasi-orderings are bqo, and {\em a fortiori\/} wqo. One 
of the most important results in this area is Laver's theorem 
(\cite{Laver}) establishing Fra\"\i ss\'e's conjecture that the 
class of countable linear orders is wqo under embeddability. For 
surveys on bqo (and wqo) theory see \cite{Milner} and 
\cite{Pou-sur}. 

We study bqos from the point of view of {\em reverse 
mathematics\/} (see \cite{SiZ2} for an overview of this research 
program; Simpson's forthcoming monograph \cite{SoSOA} is the most 
complete reference in this area) by trying to establish which 
subsystems of second order arithmetic are necessary to prove the 
theorems of bqo theory. Some results of this kind have already 
been obtained for wqo theory: see \cite{SiKr}, \cite{FRS}, 
\cite{SiHi} and \cite{Clote}. We concentrate on the very first 
result about bqos proved by Nash-Williams (\cite{NW2}) which we 
denote by \NWT\ and which asserts that if a quasi-ordering $Q$ is 
bqo then $\Qt$ (the set of all countable sequences of elements of 
$Q$ ordered by embeddability) is also bqo. Actually the very 
concept of bqo was introduced by Nash-Williams in his search for 
a notion closed under the operation $Q \mapsto \Qt$. 

Nash-Williams' original proof of \NWT\ and all its refinements 
rely on a result that in contemporary literature is called the 
{\em minimal bad array lemma}, a powerful technique that is 
widely used to prove that a quasi-ordering is a bqo: for 
treatments in which this technique is explicitly isolated see 
\cite{Lavmba}, \cite{Milner} (where is called the {\em 
forerunning technique\/}) and \cite{Sibqo}, \cite{vEMS} (which 
use a different but equivalent definition of bqo). The proof of the 
minimal bad array lemma appears to require rather strong 
set-existence axioms (e.g.\ the comprehension scheme 
characterizing the subsystem of second order arithmetic \PPCA) and 
this fact hinted that the results of bqo theory could be 
proof-theoretically very strong. In this paper we use only a weak 
version of the minimal bad array lemma which we call the {\em 
locally minimal bad array lemma\/} and which is a direct 
generalization of the minimal bad sequence lemma of wqo theory: 
this is actually the first version of the minimal bad array lemma 
proved for a specific quasi-ordering by Nash-Williams in 
\cite{NW1}. Using the locally minimal bad array lemma we prove 
\NWT\ in the system \PCA, thus showing that this theorem 
of bqo theory is not stronger than many other mathematical 
results which are known to be equivalent to \PCA\ (\cite{SiZ2}, 
\cite{SoSOA}). This result can be considered as a first step 
towards the understanding of which set-existence axioms are 
needed to prove other results of bqo theory, such as Fra\"\i 
ss\'e's conjecture, that have been object of much interest in 
recent times (\cite{Clote}, \cite{Shore}). In \cite{Clote} Clote 
claimed to have proved \NWT\ in \ATR, but unfortunately that 
proof is incorrect, as Clote himself has acknowledged (personal 
communication). 

As in \cite{thesis} (see also \cite{fbqo}) we use the 
locally minimal bad array lemma to prove the generalization of 
the classical Higman's theorem (\cite{Hig}) which states that $Q$ 
bqo implies $\Qom$ bqo and is denoted by \GHT. This is the only
step of our proof of \NWT\ where we use some form of the minimal
bad array lemma. The other main 
ingredient of our proof is hereditarily indecomposable sequences 
which are analogous to the hereditarily additively indecomposable 
countable linear orderings introduced by Laver in \cite{Laver} 
and used by Rosenstein (\cite{Ros}) in his exposition of the 
proof of Fra\"\i ss\'e's conjecture. 

In addition to proving \NWT\ in \PCA\ we show that \NWT\ is not 
equivalent to \PCA\ and we obtain some partial results toward the 
proof of \NWT\ in the weaker system \ATR. A proof of \NWT\ in 
\ATR\ would establish the logical strength of this theorem 
because it follows from the results obtained by Shore in 
\cite{Shore} that \NWT\ implies \ATR\ over the weak base system 
\RCA. The only step of the proof of \NWT\ we cannot carry out in 
\ATR\ is \GHT. We have therefore reduced the problem of proving 
\NWT\ in \ATR\ to the problem of proving \GHT\ in \ATR: since 
\GHT\ is the particular case of \NWT\ in which we deal only with 
finite sequences we believe that this is a significant step 
forward. 

We also prove that the locally minimal bad array lemma is 
equivalent (over \RCA) to \PCA\ and hence not provable in \ATR: 
therefore a proof of \GHT\ in \ATR\ (if such a proof exists) 
must avoid the use of the locally minimal bad array lemma. 

The plan of the paper is as follows. In section~\ref{sub} we 
briefly present the subsystems of second order arithmetic we will 
be dealing with and quote the results we will be using in the 
following sections: this part is not self-contained and assumes 
some acquaintance with the basic ideas of reverse mathematics. In 
section~\ref{bqo} we formally define the notion of bqo within the 
appropriate subsystems of second order arithmetic and prove in 
\PCA\ the locally minimal bad array lemma: this section is 
intended to be self-contained and we do not assume any previous 
knowledge of bqo theory. In section~\ref{NWT} we state precisely 
\NWT\ and \GHT. In section~\ref{PCA} we obtain the main result of 
the paper by showing that \NWT\ is provable in \PCA, while 
in section~\ref{ATR} we show that \NWT\ does not imply \PCA\ by 
studying the possibility of proving \NWT\ in \ATR. In 
section~\ref{rev} we prove that the locally minimal bad array 
lemma (even in its simplest form, i.e.\ the minimal bad sequence 
lemma) implies \PCA. 

Most of the results of this paper were originally obtained in the 
author's Ph.d.\ thesis (\cite{thesis}) as a consequence of a fine 
analysis of the notion of bqo. I am very much indebted to my 
advisor Stephen G. Simpson, who introduced me to bqo theory and 
to the techniques of reverse mathematics, supported me throughout 
this research and suggested several improvements to the paper: 
the results of section~\ref{rev} were indeed obtained jointly 
with him. 

\section{Subsystems of second order arithmetic}\label{sub}

The subsystems of second order arithmetic we will be using are the following.

\begin{definition}
The system \RCA\ consists of the induction scheme restricted to \S01 formulas
and the comprehension scheme restricted to \D01 formulas.

The systems \ACA, \PCA\ and \PPCA\ (in order of increasing strength) are
obtained by restricting the comprehension scheme to arithmetical, \P11 and
\P12 formulas respectively.

The system \ATR\ (which is intermediate between \ACA\ and \PCA) is obtained
by adding to \ACA\ a scheme asserting that we can iterate arithmetical
comprehension along any well-ordering.

The system \SAC\ is obtained by adding to \ACA\ the following \S11 choice
scheme:
$$\forall n\; \exists X\; \varphi(n,X) \imply \exists Y\; \forall n\;
\varphi(n,(Y)_n)$$
where $\varphi$ is a \S11 formula and $(Y)_n = \set{m}{(m,n) \in Y}$ (here
$(m,n)$ is a fixed pairing function).
\end{definition}

All these subsystems of second order arithmetic have been widely studied in
the program of {\em reverse mathematics\/} begun by Harvey 
Friedman and Stephen Simpson: more information about them can be 
found in \cite{SiZ2} and, with greater detail, in Simpson's 
forthcoming monograph (\cite{SoSOA}). Here we limit ourselves to 
the statement of a few theorems that will be used in the next 
sections. Whenever we begin a definition, lemma or theorem by the 
name of one of these subsystems between parenthesis we mean that the 
definition is given, or the statement provable, within that subsystem. 

\begin{theorem}\label{DCA}
\ATR\ proves the \D11 comprehension scheme.
\end{theorem}
\begin{pf}
See \cite{SoSOA}, section V.5.
\end{pf}

\begin{theorem}\label{SAC}
\ATR\ proves the \S11 choice scheme and therefore \ATR\ is stronger than
\SAC.
\end{theorem}
\begin{pf}
See \cite{SoSOA}, section V.8.
\end{pf}

We will need to discuss countable ordinals, which can be coded in 
the language of second order arithmetic by means of
well-orderings: details of this coding can be found in \cite{FH}, 
\cite{Hirst} and \cite{SoSOA}. In particular well-orderings are 
compared by embeddability. As usual, we will use letters from the 
beginning of the Greek alphabet to denote ordinals. In 
\cite{Hirst} are given precise definitions of all basic ordinal 
operations, including exponentiation, and the following theorem 
is proved. 

\begin{theorem}\label{Cantor}
(\ATR) Cantor's normal form theorem, i.e.\ the statement ``for 
any ordinal $\alpha$ there exists a finite sequence $\alpha_1 
\geq \dots \geq \alpha_n$ of ordinals such that $\alpha = 
\omega^{\alpha_1} + \dots + \omega^{\alpha_n}$''. 
\end{theorem}
\begin{pf}
See \cite{Hirst}, where it is shown that, over \RCA, \ATR\ is 
equivalent to Cantor's normal form theorem.
\end{pf}

\begin{definition}
(\RCA) Ordinals of the form $\omega^\alpha$ are called indecomposable.
\end{definition}

\begin{definition}
An $\omega$-model $\M$ is a model for the language of second order arithmetic
whose first order part (i.e.\ the part regarding the natural numbers) is the
standard model $\left(\N, {+}, {\cdot}, 0, 1, {<}\right)$. Therefore $\M$ can
be identified with the subset of $\cal P(\N)$ on which the set variables are
allowed to range and we write $X \in \M$ to mean that $X$ belongs to this
subset.

We write $\M \models \varphi$ to mean that the formula $\varphi$ holds within
$\M$ and we say that $\M$ is a model of the subsystem $\bold T$ of second
order arithmetic if all the axioms of $\bold T$ hold in $\M$.
\end{definition}

\begin{definition}
A $\beta$-model $\M$ is an $\omega$-model such that for any \S11 
sentence $\varphi$ with parameters from $\M$, $\M \models 
\varphi$ if and only if $\varphi$ is true.
\end{definition}

Countable $\omega$- and $\beta$-models can be coded in \RCA\ and 
several results concerning them can be proved in appropriate 
subsystems of second order arithmetic (see \cite{SoSOA}, chapters 
VII and VIII). The ones we will use in section~\ref{ATR} are the 
following. 

\begin{theorem}\label{omega}
(\ATR) For any $X \subseteq \N$ there exists a countable $\omega$-model $\M$
of \SAC\ such that $X \in \M$.
\end{theorem}
\begin{pf}
See \cite{SoSOA}, section VIII.4, or \cite{SiKo}.
\end{pf}

\begin{theorem}\label{beta}
(\PCA) For any $X \subseteq \N$ there exists a countable $\beta$-model $\M$
such that $X \in \M$.
\end{theorem}
\begin{pf}
See \cite{SoSOA}, section VII.2.
\end{pf}

The following fact is well-known but seldom explicitly stated: a slightly weaker version of it is proved in \cite{AMS}. 

\begin{corollary}\label{notP12}
No \P12 statement consistent with \ATR\ implies \PCA\ over \ATR.
\end{corollary}
\begin{pf}
Suppose \ATR\ proves that $\forall X \psi(X)$ implies \PCA, where 
$\psi$ is \S11 and suppose also that the theory $\bold T$ consisting 
of \ATR\ plus $\forall X \psi(X)$ is consistent. 

We reason in $\bold T$. $\bold T$ includes \PCA\ and by 
theorem~\ref{beta} there exists a $\beta$-model $\M$. For every 
$X \in \M$ we have $\M \models \psi(X)$: hence $\M \models 
\forall X \psi(X)$. Moreover every $\beta$-model is a model of 
\ATR\ (\cite{SoSOA}, section VII.2) and hence $\M$ is a model of 
$\bold T$. 

Thus we have shown in $\bold T$ that there exists a model of 
$\bold T$: this contradicts G\"odel's Second Incompleteness 
Theorem. 
\end{pf}

\section{Bqo theory}\label{bqo}

An alternative definition of bqo, equivalent to Nash-Williams' 
original one, has been introduced by Simpson in \cite{Sibqo}: 
this definition uses Borel maps and allows the use of many 
powerful tools from descriptive set theory in bqo theory. 
Simpson's definition has proved to be very useful but does not 
lend itself to the analysis we pursue in this paper; therefore we 
will use Nash-Williams' definition, which has a strong 
combinatorial flavor. 

\begin{definition}\label{seq}
(\RCA) If $s$ is a finite sequence of natural numbers we denote 
by $\lh(s)$ its length and, for every $i< \lh(s)$, by $s(i)$ its 
($i+1$)-th element. Then we write this sequence as $\<s(0), 
\dots, s(\lh(s) -1)\>$ or as $\<s(i) : i < \lh(s)\>$. In particular 
$\<\>$ denotes the unique sequence of length $0$.

If $s$ and $t$ are finite sequences we write $s \initeq t$ if $s$ 
is an initial segment of $t$, i.e.\ if $\lh(s) \leq \lh(t)$ and 
$\forall i < \lh(s) \; s(i)=t(i)$. We write $s \subseteq t$ if 
the range of $s$ is a subset of the range of $t$, i.e.\ if 
$\forall i < \lh(s) \; \exists j < \lh(t) \; s(i)=t(j)$. $s \init 
t$ and $s \subset t$ have the obvious meanings and we extend 
these notations also to the case where $t$ is an infinite 
sequence (i.e.\ a function from $\N$ to $\N$). 

We write $s \conc t$ for the concatenation of $s$ and $t$, i.e.\ 
the unique sequence $u$ such that $\lh(u) = \lh(s) + \lh(t)$, for 
every $i<\lh(s)$ $u(i) = s(i)$ and for every $i <\lh(t)$ 
$u(\lh(s) + i) = t(i)$. 

If $s$ is a finite sequence and $i \leq \lh(s)$ we denote by 
$s[i]$ the initial segment of $s$ of length $i$, i.e.\ the unique 
sequence $t$ such that $t \initeq s$ and $\lh(t) = i$. If $\lh(s) 
> 0$ we denote by $s^-$ the sequence obtained from $s$ by 
deleting its first element, i.e.\ $s^- = \<s(1), \dots, s(\lh(s) 
-1)\>$. These notations are extended to infinite sequences as 
well. 

If $X \subseteq \N$ we denote as usual by $\sq X$ (resp.\ $\Sq 
X$) the set of all finite (resp.\ infinite) subsets of $X$. 
Throughout this paper we will identify a subset of $\N$ with the 
unique sequence enumerating it in increasing order. Therefore if 
$s \in \sq \N$ we have that $\lh(s)$ is the cardinality of $s$ 
and, for $i< \lh(s)$, $s(i)$ is the $(i+1)$-th least element of 
$s$. Similarly if $X \in \Sq \N$ then $X(i)$ is the $(i+1)$-th 
least element of $X$.
\end{definition} 

\begin{remark}
Obviously $\Sq X$ does not formally exists within \RCA\ or
even full second order arithmetic: the use of this notation is just a
convenient shorthand. We will refer to objects like $\Sq X$ as 
classes: a class $\cal C$ is the collection of all $Y$ satisfying 
some formula $\varphi(Y)$. $ Y \in \cal C$ then stands for 
$\varphi(Y)$.
\end{remark}

\begin{definition}
(\ACA) If $B$ is a subset of $\sq \N$ we denote by $\base(B)$ the set
$$\set{n}{\exists s \in B \; \exists i < \lh(s) \; s(i) = n}$$
\end{definition}

\begin{remark}
Since \RCA\ lacks arithmetical comprehension the set $\base(B)$ does not formally exist within \RCA, but we can use this notation again as a
convenient shorthand.
\end{remark}

The following definitions are basic for bqo theory and were first given by
Nash-Williams in \cite{NW1} and \cite{NW2}, where the notion of
better-quasi-ordering was introduced.

\begin{definition}\label{blbar}
(\RCA) A set $B \subseteq \sq \N$ is a block if it satisfies the following
conditions:\begin{enumerate}
\item $\base(B)$ is infinite;
\item $\forall X \in \Sq{\base(B)} \; \exists s \in B \; s \init X$;
\item $\forall s,t \in B \; s \not\init t$.
\end{enumerate}
$B$ is a barrier if it satisfies (1), (2) and
\begin{enumerate}
\item[$(3')$] $\forall s,t \in B \; s \not\subset t$.
\end{enumerate}
\end{definition}

\begin{remark}
Since $s \init t$ implies $s \subset t$ it is immediate that in \RCA\ every
barrier is a block. Notice that if $B$ is a block and $Y \in \Sq{\base(B)}$
there exists a unique block $B' \subseteq B$ such that $\base(B') = Y$,
namely $B' = \set{s \in B}{s \subset Y}$. Notice also that if $B$ is a
barrier then $B'$ is also a barrier: we say that $B'$ is a subbarrier of
$B$.
\end{remark}

\begin{definition}
(\RCA) Let $s,s' \in \sq \N$: we write $s \tri s'$ if there exists $t \in \sq
\N$ such that $s \initeq t$ and $s' \initeq t^-$.
\end{definition}

\begin{remark}
$\tri$ is neither reflexive nor transitive: we have $\<3,5,6\> \tri
\<5,6,8,9\> \tri \<6,8\>$ but $\<3,5,6\> \ntri \<6,8\>$. In \RCA\ it is easy
to prove that if $s$ and $t$ are members of a barrier and $s \tri t$ holds
then $\lh(s) \leq \lh(t)$.
\end{remark}

We can now define the main notions of bqo theory.

\begin{definition}
(\RCA) Let $Q$ be a set (i.e.\ a subset of $\N$) and let 
$\preceq$ be a quasi-ordering on this set (i.e.\ a subset of $Q 
\times Q$ which encodes a binary relation which is reflexive and 
transitive). If $B$ is a barrier and $f: B \to Q$ we say that $f$ 
is a $Q$-array. A $Q$-array $f$ is good (with respect to $\preceq$ 
when there is danger of confusion) if there exist $s, t \in B$ 
such that $s \tri t$ and $f(s) \preceq f(t)$. If $f$ is not good 
we say that it is bad. If for every $s,t \in B$ such that $s \tri 
t$ we have $f(s) \preceq f(t)$ we say that $f$ is perfect. 
\end{definition} 

\begin{definition}\label{bqodef}
(\RCA) Let $(Q,\preceq)$ be a quasi-ordering as above: we say that $\preceq$
is a countable bqo if every $Q$-array is good with respect to $\preceq$.
\end{definition}

\begin{remark}\label{wqo}
To see that a bqo is a wqo consider the barrier $B$ of all 
singletons of $\N$. 
Then $B$ is in one-to-one correspondence with 
$\N$ and hence a $Q$-array with domain $B$ can be viewed as a 
function $f: \N \to Q$. Since $\<n\> \tri \<m\>$ if and only if 
$n<m$, $f$ is good if and only if there exist $n$ and $m$ such 
that $n<m$ and $f(n) \preceq f(m)$.
\end{remark}

\begin{remark}
We have given the definition of bqo only for countable quasi-orderings 
(but in section~\ref{NWT} we will also explain what we mean by 
saying that a class is bqo): this is due to the limitations of 
the language of second-order arithmetic, but in the case of bqos 
the restriction to the countable case does not affect the 
development of the theory. In fact, since the range of a 
$Q$-array is always countable, it is immediate that a 
quasi-ordering $Q$ is bqo if and only if every countable subset 
of $Q$ is bqo.
\end{remark}

We now show that some properties of blocks, barriers and
$Q$-arrays which were proved by Nash-Williams in \cite{NW0} and 
\cite{NW1} are actually provable within \ATR. 

\begin{theorem}\label{part}
(\ATR) If $B$ is a block (resp.\ barrier) and $B= B_1 \cup B_2$ 
then there exists a block (resp.\ barrier) $B'$ such that either 
$B' \subseteq B_1$ or $B' \subseteq B_2$.
\end{theorem}
\begin{pf}
The case in which $B$ is a block is essentially the clopen Ramsey 
theorem that can be rephrased for formalization within \ATR\ as 
``if $B_1, B_2 \subseteq \sq X$ are such that for all $Y \in \Sq 
X$ there exists $s \in B_1 \cup B_2$ such that $s \init Y$ then 
there exist $H \in \Sq X$ and $i \in \{1,2\}$ such that for all 
$Y \in \Sq H$ there exists $s \in B_i$ such that $s \init Y$''.
For proving our theorem then it suffices to let $B' \subseteq B$
be the unique block with $\base(B') = H$.
In \cite{FMS} and \cite{SoSOA} it is shown that, over \ACA, \ATR\ 
is equivalent to the clopen Ramsey theorem. 

The case in which $B$ is a barrier follows, because any block 
contained in a barrier is a barrier. Both statements are studied 
in great detail in \cite{ClRam}. 
\end{pf}

\begin{lemma}\label{blo-bar}
(\ATR) Every block contains a barrier.
\end{lemma}
\begin{pf}
We reason in \ATR. Let $B$ be a block: define $B_1 = \set{s \in 
B}{\forall t \in B\; t \not\subset s}$ and $B_2 = B \setminus 
B_1$. By theorem~\ref{part} there exists a block $B'$ such that 
either $B' \subseteq B_1$ or $B' \subseteq B_2$. In the first 
case $B'$ is a barrier and thus it suffices to show that $B' 
\subseteq B_2$ leads to a contradiction. 

If $B' \subseteq B_2$ let $s_0 \in B'$ be such that $\lh(s_0) = \min
\set{\lh(s)}{s \in B'}$. Since $s_0 \notin B_1$ there exists $t \in B$ such
that $t \subset s$. Obviously $\lh(t) < \lh(s_0)$. Since $B' \subseteq B$ is
a block we must have $t \in B'$, thereby contradicting the minimality of
$\lh(s_0)$.
\end{pf}

\begin{lemma}\label{bad/perfect}
(\ATR) Let $f:B \to Q$ be a $Q$-array. There exists a barrier $B' \subseteq
B$ such that $f \restriction B'$ is either bad or perfect.
\end{lemma}
\begin{pf}
We reason in \ATR. Let $B^2 = \set{s \cup t}{s,t \in B \land s 
\tri t}$. It is easy to check that $B^2$ is a barrier and that 
for every $t \in B^2$ there are unique $\pz(t),\po(t) \in B$ 
such that $t= \pz(t) \cup \po(t)$ and $\pz(t) \tri \po(t)$.

Let $B_1 = \set{t \in B^2}{f(\pz(t)) \preceq f(\po(t))}$ and $B_2 = B^2
\setminus B_1$. By theorem~\ref{part} there exists a barrier $B'' 
\subseteq B^2$ such that either $B'' \subseteq B_1$ or $B'' 
\subseteq B_2$. Let $B' = \set{s \in B}{s \subset \base(B'')}$: 
$B'$ is a barrier. If $B'' \subseteq B_1$ then $f \restriction 
B'$ is perfect, if $B'' \subseteq B_2$ then $f \restriction B'$ 
is bad. 
\end{pf}

To state the locally minimal bad array lemma we need the 
following definition. 

\begin{definition}\label{compatible}
(\RCA) Let $(Q,\preceq)$ be a quasi-ordering. A transitive binary 
relation $<'$ on $Q$ is compatible with $\preceq$ if for every 
$q_0, q_1 \in Q$ we have that $q_0 <' q_1$ 
implies $q_0 \preceq q_1$. We write $q_0 \leq' q_1$ for $q_0 <' 
q_1 \lor q_0 = q_1$. In this situation if $f$ and $g$ are
$Q$-arrays we write $f \leq' g$ if $\dom(f) \subseteq \dom(g)$ and 
$\forall s \in \dom(f)\; f(s) \leq' g(s)$. We write $f <' g$ if 
$f \leq' g$ and $f \neq g \restriction \dom (f)$ (i.e.\ $\exists 
s \in \dom(f)\; f(s) <' g(s)$). A $Q$-array $f$ is locally 
minimal bad with respect to $<'$ if $f$ is bad (with respect to 
$\preceq$) and there is no $g <' f$ which is bad (with respect to 
$\preceq$).
\end{definition} 

The following theorem is the locally minimal bad array lemma: the 
proof we give here is similar to Nash-Williams' original proof 
(\cite{NW1}) and rather different from the one given in 
\cite{thesis} and \cite{fbqo}: it has the advantage of being 
formalizable within \PCA\ without the use of $\beta$-models. 

\begin{theorem}\label{lmba}
(\PCA) Let $(Q, \preceq)$ be a quasi-ordering, $<'$ a well-founded 
relation on $Q$ which is compatible with $\preceq$ and $f$ a bad 
$Q$-array: then there exists $g \leq' f$ which is locally minimal 
bad with respect to $<'$.
\end{theorem}
\begin{pf}
We reason within \PCA. Given $f$ let $B = \dom(f)$ and let $T$ be the tree 
of the sequences of the form $\<(s_0,q_0),\dots,(s_{k-1},q_{k-1})\>$ 
satisfying the following conditions:
\begin{enumerate}
\item $\forall i < k \left(s_i \in B \land q_i \in Q \land \forall j < 
i\; s_j \neq s_i\right)$;
\item there exists a bad $Q$-array $g$ such that $\forall i < k 
\left(s_i \in \dom(g) \land g(s_i) = q_i\right)$ and $g \leq' f$.
\end{enumerate}
$T$ exists by \S11 comprehension (which is equivalent to \P11 
comprehension and hence available in \PCA): in fact condition (2)
in the above definition asserts the existence of a set $X \in 
\Sq{\base(B)}$ and a map $g: B \cap \sq X \to Q$ satisfying 
$$\forall s,t \in B \cap \sq X \left(g(s) \leq' f(s) \land (s 
\tri t \imply f(s) \npreceq f(t))\right)$$
together with $\forall i < k \left(s_i \subset X \land g(s_i) = 
q_i\right)$. 

We will now define recursively two infinite sequences $\<s_k\>$ 
and $\<q_k\>$ such that for every $k$ the finite sequence $t_k = 
\<(s_i,q_i): i < k \>$ belongs to $T$. Suppose that $s_0, \dots, 
s_{k-1}$ and $q_0, \dots, q_{k-1}$ have already been defined and 
notice that there exists at least one pair $(s,q) \in B \times Q$ 
such that $t_k \conc \<(s,q)\> \in T$. Let $s_k$ be such that 
$\max(s_k)$ is minimal among the $s \in B$ such that for some $q$ 
we have $t_k \conc \<(s,q)\> \in T$. Notice that if $i<k$ then 
$\max(s_i) \leq \max(s_k)$. Now let $q_k$ be minimal with respect 
to $<'$ among the $q \in Q$ such that $t_k \conc \<(s_k,q)\> \in 
T$ ($q_k$ exists because $<'$ is well-founded). 

Let $\Bt = \set{s_k}{k \in \N}$ and let $g: \Bt \to Q$ be defined 
by $g(s_k) = q_k$. We claim that $\Bt$ is a subbarrier of $B$ and $g \leq' 
f$ is a locally minimal bad $Q$-array.

To prove that $\Bt$ is a barrier we begin by noticing that 
$\base(\Bt)$ is infinite because the $s_k$'s are all distinct and 
that condition $(3')$ in the definition of barrier is obviously 
satisfied because $\Bt \subseteq B$ and $B$ is a barrier. To 
prove condition (2) in the definition of barrier we will prove 
that if $s \in B$ is such that $s \subset \base(\Bt)$ then $s \in 
\Bt$: since condition (2) holds for $B$ this suffices to show that $\Bt$ is 
a subbarrier. Suppose that for some $s \in B$ such that
$s \subset \base(\Bt)$ we have $s \notin \Bt$: since $\base(\Bt) = 
\bigcup_k s_k$ for some $k$ we have $\max(s_k) > \max (s)$ and 
$s \subseteq \bigcup_{i<k} s_i$. Let $g_k$ be a $Q$-array 
satisfying condition (2) in the definition of $T$ for the 
sequence $t_k$: since $\dom(g_k)$ is a subbarrier of $B$ we have that $s
\in \dom(g_k)$. But then $t_k \conc \<(s, g_k(s)\> \in T$ 
and this violates the minimality of $\max (s_k)$. This completes 
the proof that $\Bt$ is a barrier.

$g \leq' f$ follows immediately from the definitions of $T$ and 
$g$.

To prove that $g$ is bad let $s_k \tri s_h$: then $\max(s_k) < 
\max(s_h)$ and hence $k<h$. Let $g_{h+1}$ be a $Q$-array 
satisfying condition (2) in the definition of $T$ for the 
sequence $t_{h+1}$: we have that $s_k, s_h \in \dom(g_{h+1})$ and, 
since $g_{h+1}$ is bad, $q_k \npreceq q_h$, i.e.\ $g(s_k) \npreceq 
g(s_h)$.

To prove that $g$ is locally minimal bad suppose that $g' <' g$ 
is bad and let $B' = \dom(g') \subseteq \Bt$. Let $k$ be minimal 
such that $g'(s_k) <' q_k$, so that for all $i < k$ such that 
$s_i \in B'$ we have $g'(s_i) = q_i$. Let $B''$ be the unique 
subbarrier of $\Bt$ such that $\base(B'') = \base (B') \cup 
\set{n \in \base (\Bt)}{n \leq \max (s_k)}$ and let $A = \set{s 
\in B''}{s \in B' \land \max(s) \geq \max (s_k)}$. Define $g'': 
B'' \to Q$ by 
$$g''(s) = \begin{cases}
g'(s) &\text{if $s \in A$}\\
g (s) &\text{if $s \notin A$} \end{cases}$$
We claim that $g''$ is bad. To see this suppose $s,t \in B''$ are 
such that $s \tri t$. There are three different possibilities:\begin{enumerate}
\item If $s \in A$ then also $t \in A$ (because $B'$ is a barrier 
and $\max(s) < \max (t)$) and $g'(s) \npreceq g'(t)$ implies 
$g''(s) \npreceq g''(t)$. 
\item If $s \notin A$ and $t \notin A$ then $g(s) \npreceq g(t)$
implies $g''(s) \npreceq g''(t)$.
\item If $s \notin A$ and $t \in A$ then $g'(t) \leq' g(t)$ and hence
$g'(t) \preceq g(t)$. On the other hand $g(s) \npreceq g(t)$ and hence
$g(s) \npreceq g'(t)$, that is $g''(s) \npreceq g''(t)$.
\end{enumerate}
Therefore $g''$ is bad. If $i \leq k$ we have $s_i 
\subseteq \set{n \in \base (\Bt)}{n \leq \max (s_k)}$ and hence $s_i \in 
B''$. Moreover $\forall i < k\; g''(s_i) = q_i$ and thus
$t_k \conc \<(s_k, g''(s_k))\> \in T$. Then $g''(s_k) = 
g'(s_k) <' q_k$ contradicts the minimality of $q_k$. This 
completes the proof that $g$ is locally minimal bad.
\end{pf}

In section~\ref{rev} we will reverse this theorem by showing that 
the locally minimal bad array lemma implies \PCA\ and hence 
cannot be proved in any weaker system.

\section{Transfinite sequences}\label{NWT}

In this section we give the precise statements of the two theorems
we will study in the next sections.

\begin{definition}
(\RCA) Given a set $Q$ we denote by $\Qt$ the class of all 
countable transfinite sequences of elements of $Q$, i.e.\ 
functions with domain a countable ordinal and range a subset of 
$Q$. $\Qom$ consists of all elements of $\Qt$ whose domain is a 
finite ordinal: these finite sequences can be viewed as finite 
sequences of natural numbers and hence $\Qom$ is actually a set 
in \RCA. 
\end{definition}

Since the property of being an ordinal is \P11 the class $\Qt$ is 
defined by a \P11 formula. 

We extend the notations of definition~\ref{seq} also to elements 
of $\Qt$ and in particular we write $\lh(\s)$ in place of 
$\dom(\s)$ whenever $\s \in \Qt$. Operations such as 
concatenation of elements of $\Qt$ require operations such as sum 
of ordinals, while notions such as $\s \initeq \tau$ are defined 
modulo comparability of ordinals. 

\begin{definition}
(\RCA) If $Q$ is equipped with a quasi-ordering $\preceq$ and 
$\s, \tau \in \Qt$ we say that $\s$ is embeddable in $\tau$ and 
write $\s \leq \tau$ if there exists a strictly increasing map 
$f: \lh(\s) \to \lh(\tau)$ such that for every $\alpha < 
\dom(\s)$ we have $\s(\alpha) \preceq \tau(f(\alpha))$. $f$ is 
called an embedding of $\s$ in $\tau$.
\end{definition} 

\begin{remark}
It is easy to see (even within \RCA) that if $\preceq$ is a 
quasi-ordering on $Q$ than $\leq$ is a quasi-ordering on $\Qt$. 
In \RCA\ the formula $\s \leq \tau$ is \S11. 
\end{remark}

\begin{lemma}\label{D11}
(\ATR) If $\s, \tau \in \Qt$ the formula $\s \leq \tau$ is \D11.
\end{lemma}
\begin{pf}
To prove that $\s \leq \tau$ is \P11 in \ATR\ we construct
a canonical (possibly partial) embedding of $\s$ in 
$\tau$. Define $f:\lh(\s) \to \lh(\tau)$ by transfinite recursion 
letting, for each $\alpha < \lh(\s)$, $f(\alpha)$ to be the least 
$\beta < \lh(\tau)$ such that $\beta > f(\gamma)$ for all $\gamma 
< \alpha$ and $\s(\alpha) \preceq \tau(\beta)$ if such a $\beta$ 
exists; $f(\alpha)$ is undefined otherwise. Let $\varphi 
(f,\s,\tau)$ be the arithmetical formula asserting that $f$ is 
built this way: then $\s \leq \tau$ is equivalent to $\forall f 
(\varphi(f,\s,\tau) \imply \dom(f) = \lh(\s))$, which is \P11. 
\end{pf} 

One of the first results in wqo theory (rediscovered many times, 
see \cite{Kruskal}) is Higman's theorem (\cite{Hig}): if $Q$ is 
wqo then $\Qom$ is wqo. From the point of view of reverse 
mathematics it follows from the results in section~4 of  
\cite{SiHi}, that Higman's theorem is provable in \ACA\ (see 
\cite{Clote} and \cite{thesis} for details) and in \cite{Clote} 
it is proved that Higman's theorem implies \ACA. In this paper we 
will be interested in the following natural generalization of 
Higman's theorem. 

\begin{definition}
(\RCA) We denote by \GHT\ the following statement: ``if $Q$ is a 
countable bqo then $\Qom$ is a countable bqo''.
\end{definition}

\GHT\ is easily provable by a generalization of the standard argument for
Higman's theorem, using some form of the minimal bad array lemma (we will do
this in lemma~\ref{GHT}).

$\Qt$ is not a set, so we cannot apply the definition of 
countable bqo to it. We can however imitate
definition~\ref{bqodef} as follows.

\begin{definition}
(\RCA) A $\Qt$-array is a sequence $\<f(s) : s \in B\>$ of 
elements of $\Qt$ indexed by a barrier $B$; a $\Qt$-array is good 
if $\exists s,t \in B (s \tri t \land f(s) \leq f(t))$. If every 
$\Qt$-array is good we say that $\Qt$ is bqo. 
\end{definition} 
                 
The next definition states the main object of interest of the
paper: Nash-Williams' theorem.

\begin{definition}
(\RCA) We denote by \NWT\ the following statement: ``if $Q$ is a 
countable bqo then $\Qt$ is bqo''.
\end{definition}

\NWT\ was proved by Nash-Williams in \cite{NW2} and
can be considered the first result of bqo theory. Actually Nash-Williams
considered arbitrary transfinite sequences of elements of $Q$, while in $\Qt$
we have only sequences of countable length: this is necessary to formalize
\NWT\ in second order arithmetic, where we cannot discuss uncountable
ordinals. Pouzet (\cite{Pou}) proved that $Q$ bqo is actually 
equivalent to $\Qt$ wqo, where $\Qt$ can be taken either in
Nash-Williams' original meaning or in our countable meaning: this 
shows that our restriction is not very relevant. 

To prove \NWT\ we will use the class of hereditarily 
indecomposable sequences and to define it we need the notion
of quasi-monotonic sequence.

\begin{definition}
(\RCA) An infinite sequence $\<\s_n : n \in \N\>$ of elements of $\Qt$ 
is quasi-monotonic if $\forall n \; \exists m>n \; \s_n \leq 
\s_m$. We denote by $\sum_n \s_n$ the infinite concatenation 
of the $\s_n$'s, i.e.\ $\sum_n \s_n = \s_0 \conc \s_1 \conc 
\dots \conc \s_n \conc \dots$.
\end{definition}

\begin{lemma}\label{notemb}
(\ATR) Suppose $\s, \tau \in \Qt$ are such that $\s = \sum_n \s_n$, $\tau =
\sum_n \tau_n$, the sequence $\<\tau_n\>$ is quasi-monotonic and $\s \nleq 
\tau$. Then there exists $n_0$ such that for all $m$ we have 
$\s_{n_0} \nleq \sum_{n \leq m} \tau_n$.
\end{lemma}
\begin{pf}
We reason in \ATR. Suppose that for every $n_0$ there is $m$ such 
that $\s_{n_0} \leq \sum_{n\leq m} \tau_n$: using the
quasi-monotonocity of $\<\tau_n\>$ we have that for every $n_0$ 
and $i$ there exists $m$ such that $\s_{n_0} \leq \sum_{i< n\leq 
m} \tau_n$. By lemma~\ref{D11} (and using the fact that in \ATR\ 
we have \D11 comprehension by theorem~\ref{DCA}) we can find the 
least such $m$. Using the \S11 choice scheme (available in \ATR\ 
by theorem~\ref{SAC}) we can find a sequence of embeddings that 
for every $n_0$ and $i$ witness $\s_{n_0} \leq \sum_{i< n\leq m} 
\tau_n$. By appropriately connecting some of these embeddings we 
construct an embedding that shows $\s \leq \tau$, a 
contradiction.
\end{pf} 

The class of hereditarily indecomposable sequences is denoted by 
$\H$ and its intuitive definition is the following: all sequences 
of length $1$ are in $\H$ and if $\<\s_n\>_{n<\omega}$ is a 
quasi-monotonic sequence of elements of $\H$ then $\sum_n \s_n$ 
is also in $\H$. This class is formalized within second order 
arithmetic by means of well-founded trees which encode the 
construction of the elements of $\H$ from sequences of length $1$ 
by repeated quasi-monotonic sums. In definition~\ref{H} we will 
not explicitly state that the tree is well-founded because, as 
noted in remark~\ref{HS11}, this follows from the other conditions of 
the definition. 

The following definition introduces some notation that is useful 
for the formalization of $\H$ within \RCA\ and which will be used 
also in sections~\ref{PCA} and~\ref{ATR}. 

\begin{definition}
(\RCA) If $\s \in \Qt$ and $\alpha < \beta \leq \lh(\s)$ let
$$\Int\alpha\beta = \set{\gamma < \lh(\s)}{\alpha \leq \gamma < 
\beta}$$
$\Int\alpha\beta$ is well-ordered by $\leq$ and hence can be 
regarded as an ordinal. By $\ss\alpha\beta$ we will denote the 
restriction of $\s$ to this ordinal.
\end{definition}

\begin{definition}\label{H}
(\RCA) $\s \in \Qt$ is hereditarily indecomposable if there 
exist a tree $T$ and two sequences $\<\alpha_t: t \in T\>$ and 
$\<\beta_t: t \in T\>$ such that the following conditions are satisfied:
\begin{enumerate}
\item $\forall t \in T\; \alpha_t < \beta_t \leq \lh(\s)$.
\item $\alpha_{\<\>} = 0$ and $\beta_{\<\>} = \lh (\s)$.               
\item If $t \in T$ is not an endnode of $T$  (i.e.\ $\exists n\; 
t \conc \<n\> \in T$) then for every $n$ we have $t \conc \<n\> 
\in T$ and $\beta_{t \conc \< n \>} = \alpha_{t \conc \< n+1 
\>}$, $\< \ss{\alpha_{t \conc \<n\>}}{\beta_{t \conc \<n\>}} 
: n \in \N\>$ is quasi-monotonic, $\alpha_t = \alpha_{t\conc 
\<0\>}$ and $\beta_t = \sup_n \beta_{t\conc \<n\>}$.
\item If $t \in T$ is an endnode of $T$ then $\beta_t = \alpha_t 
+ 1$. 
\end{enumerate}
The class of all hereditarily indecomposable sequences is denoted
by $\H$.
\end{definition}

\begin{remark}\label{HS11}
Notice that if $T$ is the tree in the above definition we can 
prove within \RCA, using the fact that $\lh(\s)$ is an ordinal, 
that $T$ is well-founded. Notice also that if $\s \in \H$ then in 
\ATR\ we can prove that $\lh(\s)$ is indecomposable. Moreover in 
\SAC\ we have that if $\s \in \Qt$ then the formula $\s \in \H$ 
is \S11 (i.e.\ $\H$ is \S11 on \P11).
\end{remark} 

Hereditarily indecomposable sequences are analogous to the hereditarily
additively indecomposable countable linear orderings introduced by Laver
(\cite{Laver}) and used by Rosenstein (\cite{Ros}, section~10.5) in his
exposition of the proof of Fra\"\i ss\'e's conjecture. Indeed the proofs of
lemmas~\ref{concH} and~\ref{HQt} (but not that of lemma~\ref{QomH}) are just
translations of the corresponding proofs given by Rosenstein for linear
orderings.

\section{\PCA\ proves \NWT}\label{PCA}

We will prove theorem~\ref{main}, the main results of the paper, by
establishing a chain of implications. Throughout the remainder of 
the paper $Q$ will denote a set equipped with a quasi-ordering $\preceq$. 

\begin{lemma}\label{GHT}
\PCA\ proves \GHT.
\end{lemma}
\begin{pf}
We reason within \PCA\ and suppose $Q$ is a countable bqo such 
that $\Qom$ is not bqo so that there exists a bad $\Qom$-array. 
On $\Qom$ $\init$ is well-founded and compatible with $\leq$. 
Hence, by theorem~\ref{lmba}, let $f: B \to \Qom$ be locally 
minimal bad with respect to $\init$. 

Clearly $\forall s \in B\; f(s) \neq \< \>$ and we can define maps
$f':B \to \Qom$ and $g:B \to Q$ such that for every $s \in B$ we have $f(s)
= f'(s) \sconc \<g(s)\>$. Since $Q$ is bqo lemma~\ref{bad/perfect} implies
that there exists a barrier $B' \subseteq B$ such that $g \restriction B'$
is perfect. The
badness of $f$ then requires $f' \restriction B'$ to be bad, contradicting
the minimality of $f$ since clearly $f' \restriction B' \init f$.
\end{pf}

The major difference between our proof and the usual proofs of 
\NWT\ lies in the way we establish the following lemma. Our proof
(which is a formalization of the proof appearing in \cite{thesis} 
and \cite{fbqo}) avoids the use of any form of the minimal bad array lemma 
and can be carried out in \ATR: in some sense it is more constructive
than the usual ones. 

\begin{lemma}\label{QomH}
(\ATR) If $\Qom$ is a countable bqo then $\H$ is bqo.
\end{lemma}
\begin{pf}
We reason within \ATR\ and suppose that $\H$ is not bqo: let $f: 
B \to \H$ be a bad $\H$-array, i.e.\ a $\Qt$-array such that for 
every $r \in B$ we have $f(r) \in \H$. We can assume that 
$\base(B) = \N$ because every barrier is isomorphic to a barrier 
with base $\N$.

According to definition~\ref{H} for every $r \in B$ there exists 
a tree and two sequences of ordinals which encode the 
construction of $f(r)$ from sequences of length $1$ by 
quasi-monotonic sums. By \S11 choice there exists sets coding a 
sequence of trees $\<T_r : r \in B\>$ and two sequences of 
sequences of ordinals $\<\<\alpha^r_t : t \in T_r\> : r \in B\>$ 
and $\<\<\beta^r_t : t \in T_r\> : r \in B\>$ which encode all 
this information. Let 
\begin{align*}
H &= \set{(r,t)}{r \in B \land t \in T_r}\\
H_0 &= \set{(r,t) \in H}{t \text{ is not an endnode of }T_r}\\
K &= \set{(r,t,i)}{r \in B \land t \in T_r \land i \geq 1 \land \alpha^r_t
+ i \leq \lh(f(r))}\end{align*}
If $r \in B$ and $f(r) = \s$ we will view $(r,t) \in H$ as a code 
for the sequence $\ss{\alpha^r_t}{\beta^r_t} \in \H$ and 
$(r,t,i) \in K$ as a code for the sequence 
$\ss{\alpha^r_t}{\alpha^r_t+i} \in \Qom$. With this coding and using \D11 
comprehension we can define a quasi-ordering $\leq^*$ on $H \cup 
K$ that corresponds to $\leq$ on the sequences coded by the elements of $H 
\cup K$. Notice that if $(r,t) \in H \setminus H_0$ then 
$(r,t)$ and $(r,t,1) \in K$ code the same sequence and 
hence are equivalent with respect to $\leq^*$. We further extend our 
identification and if $(r,t) \in H_0$ we write $(r,t) = \sum_n 
(r,t\conc\<n\>)$: obviously this sum is quasi-monotonic with 
respect to $\leq^*$. 

We will define by arithmetic recursion a sequence of functions 
$f_k: B_k \to H_0 \cup K$. We will then prove by induction that 
$B_k$ is a block with $\base(B_k) = \N$ and $f_k$ is bad (the 
definition of bad makes sense also if $B_k$ is a block instead of 
a barrier) with respect to $\leq^*$. 

We start with $B_0 = B$ and for every $r \in B$ we define $f_0(r) 
= (r,\<\>,1) \in K$ if $(r,\<\>) \in H \setminus H_0$, $f_0(r) = 
(r,\<\>) \in H_0$ otherwise. Recall the definitions of $B^2$, 
$\pz(s)$ and $\po(s)$ from the proof of lemma~\ref{bad/perfect} 
and suppose we have defined $f_k$ and $B_k$: let
$$C_k = \set{s \in B_k}{f_k(s) \in K} \quad \text{and} \quad D_k 
= \set{s \in B_k^2}{f_k(\pz(s)) \in H_0}\text{.}$$
Set $B_{k+1} = C_k \cup D_k$. 

Define $f_{k+1}: B_{k+1} \to H_0 \cup K$ by cases as follows:\begin{enumerate}
\item if $s \in C_k$ let $f_{k+1} (s) = f_k (s)$.
\item if $s \in D_k$, $f_k(\po(s)) \in H_0$ and $(r,t) = 
f_k(\pz(s)) \nleq^* f_k(\po(s)) = (r',t')$ we have $f_k(\pz(s)) = 
\sum_n (r,t\conc\<n\>)$ and $f_k(\po(s)) = \sum_n 
(r',t'\conc\<n\>)$ with both sums quasi-monotonic. By 
lemma~\ref{notemb} (and the definition of $\leq^*$) there exists 
$n_0$ minimal such that for every $m$ we have $(r,t\conc\<n_0\>) 
\nleq^* \sum_{n\leq m} (r',t'\conc\<n\>)$. If $t\conc\<n_0\>$ is 
not an endnode of $T_r$ let $f_{k+1}(s) = (r,t\conc\<n_0\>) \in 
H_0$, otherwise let $f_{k+1}(s) = (r,t\conc\<n_0\>,1) 
\in K$ (which, as noted above, is equivalent with respect to
$\leq^*$ to $(r,t\conc\<n_0\>) \in H \setminus H_0$). 
\item if $s \in D_k$, $f_k(\po(s)) \in H_0$ and $f_k(\pz(s)) 
\leq^* f_k(\po(s))$ we let $f_{k+1}(s) = f_k(\pz(s))$. It will 
follow from the badness of $f_k$ that this case never occurs. 
\item if $s \in D_k$, $f_k(\pz(s)) = (r,t)$ and $f_k(\po(s)) = 
(r',t',i) \in K$ let $f_{k+1}(s) = (r,t,i+1) \in K$ (notice that
since $(r,t) \in H_0$ we have $\alpha^r_t + \omega \leq \lh(f(r))$. 
\end{enumerate}
This completes the definitions of $B_k$ and $f_k$. We shall now 
prove by induction that they have the desired properties.

It is easy to check inductively that $\base(B_k) = \N$ and that 
the elements of $B_k$ are incomparable under $\init$. Since in 
\ATR\ we do not have \P11 induction we cannot prove inductively 
that for all $k$ $B_k$ satisfies condition (2) of 
definition~\ref{blbar}. This difficulty however can be overcome 
as follows: we fix $X \in \Sq\N$ and prove by arithmetical 
induction on $k$ that $\forall m\; \exists s \in B_k\; s \init 
X^{-m}$, where $X^{-m}$ denotes the result of iterating $X 
\mapsto X^-$ $m$ times (i.e.\ $X^{-m} = \< X(m), X(m+1), 
\ldots\>$). Supposing this holds for $B_k$, let $s \in B_k$ be 
such that $s \init X^{-m}$: if $f_k(s) \in K$ then $s \in C_k 
\subseteq B_{k+1}$; if $f_k(s) \in H_0$ let $t \in B_k$ be such 
that $t \init X^{-(m+1)}$, so that $s \tri t$, $s \cup t \in D_k 
\subseteq B_{k+1}$ and $s \cup t \init X^{-m}$. 

We now assume that $f_k$ is bad (and hence case (3) never occurs 
in the definition of $f_{k+1}$) and prove that $f_{k+1}$ is bad. 
Let $s, s' \in B_{k+1}$ be such that $s \tri s'$. We need to 
distinguish four different cases:\begin{enumerate} 
\item[(a)] if $s,s' \in C_k$ then $f_{k+1}(s) = f_k(s) \nleq^* 
f_k(s') = f_{k+1}(s')$. 
\item[(b)] if $s \in C_k$ and $s' \in D_k$ then $s \tri \pz(s')$ and 
hence $f_k(s) \nleq^* f_k(\pz(s'))$. Moreover it is clear that 
$f_{k+1}(s') \leq^* f_k(\pz(s'))$ and therefore $f_{k+1}(s) = 
f_k(s) \nleq^* f_{k+1}(s')$. 
\item[(c)] if $s \in D_k$ and $s' \in C_k$ we have $\po(s) = s'$ and 
hence in defining $f_{k+1}(s)$ we apply case (4). Therefore the 
length of the sequence coded by $f_{k+1}(s)$ is strictly greater 
than the length of the sequence coded by $f_k(s') = f_{k+1}(s')$ 
and we necessarily have $f_{k+1}(s) \nleq^* f_{k+1}(s')$. 
\item[(d)] if $s,s' \in D_k$ we have $\po(s) = \pz(s')$ and hence 
$f_{k+1}(s)$ is defined applying case (2). Therefore if 
$f_k(\pz(s')) = (r,t)$ for every $m$ we have that $f_{k+1}(s) 
\nleq^* \sum_{n\leq m} (r,t\conc\<n\>)$. On the other hand,
since $f_{k+1}(s')$ is defined using either case (2) or case 
(4), for some $m$ $f_{k+1}(s') \leq^* \sum_{n\leq m} 
(r,t\conc\<n\>)$. These two facts imply $f_{k+1}(s) \nleq^* 
f_{k+1}(s')$.
\end{enumerate} 

Notice that if $X \in \Sq\N$ and for every $k$ we define $s_k(X)$ 
to be the unique sequence such that $s_k(X) \in B_k$ and $s_k(X) 
\init X$, we have either $s_k(X) \in C_k$ or that the length of 
the sequence coded by $f_{k+1}(s_{k+1}(X))$ is strictly smaller 
than the length of the sequence coded by $f_k(s_k(X))$. 

Now let us define $B = \bigcup_k C_k$. We claim that $B$ is a 
block. To see this notice that, by the above observation, for any 
$X \in \Sq\N$ there exists $k$ such that $s_k(X) \in C_k$: this 
implies that $\base(B) = \N$ and that condition (2) in 
definition~\ref{blbar} is satisfied. Condition (3) follows from 
the fact that if $s_k(X) \in C_k$ then for all $k' > k$ we have 
$s_{k'}(X) = s_k(X)$. 

Define $\fom: B \to K$ by setting $\fom(s) = f_k(s)$ for the 
least (or any) $k$ such that $s \in C_k$. It is immediate to 
check that $\fom$ is bad and hence if $B' \subseteq B$ is the 
barrier given by lemma~\ref{blo-bar} $\fom \restriction B'$ is 
also bad. Since $\fom$ clearly codes a bad $\Qom$-array we have 
that $\Qom$ is not a countable bqo and the proof is complete. 
\end{pf} 

The following result will be improved by lemma~\ref{concATR}, but in the
present section it suffices to prove theorem~\ref{main}.

\begin{lemma}\label{concH}
(\PCA) If $\H$ is wqo then every element of $\Qt$ can be written 
as the concatenation of finitely many elements of $\H$. 
\end{lemma}
\begin{pf}
Let $\s \in \Qt$ be given: we prove by transfinite induction on 
$\Int\alpha\beta$ that the lemma holds for $\ss\alpha\beta$. For 
$\alpha = 0$, $\beta = \lh(\s)$ this establishes the lemma for 
$\s$. For $\ss\alpha\beta$ the lemma states
$$\exists \<\alpha_i: i \leq n\> \left( \alpha_0 = \alpha \land 
\alpha_n = \beta \land \forall i < n (\alpha_i < \alpha_{i+1} 
\land \ss{\alpha_{i}}{\alpha_{i+1}} \in \H) \right)$$
which by remark~\ref{HS11} is \S11. Hence this induction can be 
carried out in \PCA. 

By theorem~\ref{Cantor} Cantor's normal form theorem is provable 
in \ATR\ and hence in \PCA. Thus it suffices to prove the 
statement for $\alpha, \beta$ such that $\Int\alpha\beta$ is 
indecomposable. 

We reason within \PCA\ and begin by noticing that the result is 
obvious when $\beta = \alpha+1$. Given $\Int\alpha\beta$ 
indecomposable and greater than $1$ we can write $\Int\alpha\beta 
= \bigcup_n \Int{\alpha_n}{\beta_n}$ so that $\ss\alpha\beta = 
\sum_n \ss{\alpha_n}{\beta_n}$. Let us write $\s_n$ in place of 
$\ss{\alpha_n}{\beta_n}$. By induction hypothesis every $\s_n$ 
can be written as the concatenation of finitely many elements of 
$\H$ and hence we can assume that for every $n$ we have $\s_n \in 
\H$. 

We claim that there exists $m$ such that $\forall n>m\; \exists 
n'>n\; \s_n \leq \s_{n'}$. Indeed if this is not the case we have 
$\forall m\; \exists n>m\; \forall n'>n\; \s_n \nleq \s_{n'}$. 
Then we can find a set $A \subset \N$ such 
that $n,n' \in A$ and $n < n'$ imply $\s_n \nleq \s_{n'}$, 
against the fact that $\H$ is wqo. 

The claim we just proved shows that for some $m$ the sequence 
$\<\s_n : n>m\>$ is quasi-monotonic and hence $\sum_{n>m} \s_n 
\in \H$. Therefore $\s_0 \sconc \dots \conc \s_m \sconc 
\left(\sum_{n>m} \s_n \right)$ expresses $\ss\alpha\beta$ as the 
concatenation of finitely many elements of $\H$. 
\end{pf}

\begin{lemma}\label{HQt}
(\PCA) If $\H$ is bqo then $\Qt$ is bqo.
\end{lemma}
\begin{pf}
We reason within \PCA. If $f: B \to \Qt$ is a $\Qt$-array for 
every $s \in B$ we can write, by lemma~\ref{concH}, $f(s) = 
\tau_{s,0} \sconc \dots \conc \tau_{s,k_s}$ where $\tau_{s,i} \in 
\H$. Since $\H$ is bqo the set $R = \set{(s,i)}{s \in B \land i 
\leq k_s}$ with the quasi-ordering defined by setting $(s,i) 
\preceq_R (t,j)$ if and only if $\tau_{s,i} \leq \tau_{t,j}$ is a 
countable bqo. We can apply lemma~\ref{GHT} and $R^{<\omega}$ is 
also a countable bqo: hence there exist $s \tri t$ in $B$ such 
that $\<(s,0), \dots, (s,k_s)\> \leq \<(t,0), \dots, (t,k_t)\>$ 
in $R^{<\omega}$ and therefore $f(s) \leq f(t)$. 
\end{pf}

\begin{theorem}\label{main}
\PCA\ proves \NWT.
\end{theorem}
\begin{pf}
It suffices to concatenate lemmas~\ref{GHT}, \ref{QomH} and~\ref{HQt}.
\end{pf}

\section{\NWT\ in \ATR}\label{ATR}

The following theorem, due to Shore, sheds further light on the 
logical strength of \NWT. 

\begin{theorem}\label{Shore}
(\RCA) The statement ``the class of countable well-orderings is wqo under
embeddability'' implies \ATR.
\end{theorem}
\begin{pf}
See \cite{Shore}.
\end{pf}

\begin{corollary}\label{lb}
\RCA\ proves that \NWT\ implies \ATR.
\end{corollary}
\begin{pf}
Let $Q$ be a quasi-ordering consisting of a single element: $Q$ 
is a countable bqo (with respect to the identity relation) and by 
\NWT\ $\Qt$ is also bqo and hence wqo. Moreover in this case 
$\Qt$ is isomorphic to the class of countable well-orderings 
ordered by embeddability. Thus the statement contained in 
theorem~\ref{Shore} is a consequence of \NWT\ in \RCA\ and \NWT\ 
implies \ATR.
\end{pf}

Thus theorem~\ref{main} and corollary~\ref{lb} establish \PCA\ 
and \ATR\ as upper and lower bounds for the logical strength of 
\NWT. In this section we show that \NWT\ is not equivalent to 
\PCA\ and discuss the possibility of proving the following 
conjecture, which would improve the result of section~\ref{PCA}. 
                          
\begin{conjecture}\label{conj}
\ATR\ proves \NWT.
\end{conjecture}

In light of corollary~\ref{lb}, proving conjecture~\ref{conj} 
would imply that the answer to the question ``which set existence 
axioms are needed to prove \NWT?'' is Arithmetical Transfinite 
Recursion. Answers of this kind are typical reverse mathematics 
results. 

The following lemma improves lemma~\ref{concH} and is a step 
towards proving conjecture~\ref{conj}. The proof consists of 
carrying out some parts of the argument used to prove 
lemma~\ref{concH} inside an $\omega$-model. 

\begin{lemma}\label{concATR}
(\ATR) If $\H$ is wqo then every element of $\Qt$ can be written as the
concatenation of finitely many elements of $\H$.
\end{lemma}
\begin{pf}                 
We reason within \ATR\ and let $\s \in \Qt$: using \S11 choice we can 
construct a single set $X$ which codes $\s$, a witness for the Cantor's 
normal form for every ordinal $\Int\alpha\beta$ with $\alpha < 
\beta \leq \lh(\s)$ and the canonical witness of the proof of 
lemma~\ref{D11} for $\ss\alpha\beta \leq \ss{\alpha'}{\beta'}$ 
(i.e.\ the unique $f$ such that $\varphi (f, \ss\alpha\beta, 
\ss{\alpha'}{\beta'})$ holds) for every two such pairs $\alpha, 
\beta$ and $\alpha', \beta'$. 

By theorem~\ref{omega} there exists a countable $\omega$-model 
$\M$ of \SAC\ such that $X \in \M$. By our construction the 
formulas $\ss\alpha\beta \leq \ss{\alpha'}{\beta'}$ are \D11 
within $\M$. 

By transfinite induction on the length of $\Int\alpha\beta$ we 
prove that
\begin{multline*}
\forall \alpha,\beta < \lh(\s) (\alpha < \beta \imply \\
\M \models (\ss \alpha\beta \text{ is the concatenation of finitely many
elements of $\H$}))
\end{multline*}
This transfinite induction can be carried out in \ATR\ because 
$\M \models \psi$ is an arithmetical formula even when $\psi$ is 
not arithmetical. Moreover this induction will complete the proof 
of the lemma because for $\alpha = 0$, $\beta = \lh(\s)$ we 
obtain that $\s$ satisfies the lemma within $\M$ and, since the 
statement of the lemma is \S11 (as already noticed in the proof 
of lemma~\ref{concH}), it satisfies the lemma also in general. 

Since we have enough instances of Cantor's normal form theorem 
within $\M$, it suffices to consider the case where 
$\Int\alpha\beta$ is indecomposable. If $\Int\alpha\beta$ has 
length $1$ the statement is trivial. Otherwise we reason in
$\M$ and write $\ss \alpha\beta = \sum_n \s_n$ as in the proof 
of lemma~\ref{concH}: by induction hypothesis we may assume that 
for every $n$ $\s_n \in \H$. 

Now we claim that $\exists m\; \forall n>m\; \exists n'>n\; 
\s_n \leq \s_{n'}$. If this is not the case we have
$$\M \models \left(\forall m\; \exists n>m\; \forall n'>n\;
\s_n \nleq \s_{n'}\right)$$
Since the formula under consideration is \D11 (to establish this 
we need to use \SAC\ so that the numerical quantifiers can be 
brought inside the set quantifier of $\nleq$) it is absolute and 
holds also outside $\M$ and, as in the proof of 
lemma~\ref{concH}, we reach a contradiction with the hypothesis 
that $\H$ is wqo. 

From the claim we complete the induction step exactly as we did in the proof
of lemma~\ref{concH} and the proof is complete.
\end{pf}

\begin{theorem}\label{status}
\ATR\ proves that \GHT\ and \NWT\ are equivalent.
\end{theorem}
\begin{pf}
Since $\Qom$ is a subset of $\Qt$ it is obvious that \NWT\
implies \GHT\ even in \RCA.

We now reason within \ATR. Assuming \GHT\ and using lemma~\ref{QomH} 
we have that if $Q$ is a countable bqo then $\H$ is bqo. 
Repeating the proof of lemma~\ref{HQt}, from \GHT\ and 
lemma~\ref{concATR} it follows that $\H$ bqo implies $\Qt$ bqo. 

Therefore using \GHT\ we can prove in \ATR\ that if $Q$ is a 
countable bqo then $\Qt$ is bqo, i.e.\ \NWT.
\end{pf}

Hence we have reduced conjecture~\ref{conj} to the following 
conjecture.

\begin{conjecture}
\ATR\ proves \GHT.
\end{conjecture}

Since $\Qom$ is a ``small'' subset of $\Qt$ this is a significant 
step towards establishing conjecture~\ref{conj}. 

Note that in the proof of lemma 15 in \cite{Clote} Clote claims 
that \GHT\ is provable even in \ACA\ by modifying the techniques 
used in \cite{SiHi} to show that Higman's theorem is provable in 
\ACA. We have not been able to find such a modification. 

A natural way of refuting our conjectures would be to prove that 
\NWT\ (or, equivalently, \GHT) implies \PCA\ over \ATR. The following
theorem shows that this is not possible.

\begin{theorem}
\ATR\ does not prove that \NWT\ implies \PCA.
\end{theorem}
\begin{pf}
Assume \ATR\ proves that \NWT\ implies \PCA. By 
theorem~\ref{status} \ATR\ proves that \GHT\ implies \PCA. \GHT\
is a \P13 statement and hence we cannot apply corollary~\ref{notP12}
directly; however, going 
back to the proof of lemma~\ref{GHT} we can observe that a 
statement $\psi$ slightly stronger than \GHT\ was obtained in \PCA\
(and hence is consistent with \ATR). 
$\psi$ states: ``if $B$ is a barrier and $Q$ is $B$-wqo then 
$\Qom$ is $B$-wqo'' where $B$-wqo means that every array with 
domain a subbarrier of $B$ (including $B$ itself) is good (this 
notion is strictly connected with that of \mbox{$\alpha$-wqo} studied in 
\cite{Pou}, \cite{Clote}, \cite{thesis} and \cite{fbqo}). \RCA\ 
proves that $B$-wqo for every barrier $B$ is equivalent to bqo and
hence that $\psi$ implies \GHT. Moreover $\psi$ is \P12 and the
fact that \ATR\ proves that it implies \PCA\ contradicts 
corollary~\ref{notP12}. 
\end{pf}

\section{The locally minimal bad array lemma implies \PCA}\label{rev} 

In this section we obtain a typical reverse mathematics result by 
showing that the locally minimal bad array lemma (i.e.\ the 
statement proved within \PCA\ in theorem~\ref{lmba}) is 
equivalent to \PCA\ over \RCA. We actually prove that the special 
case obtained by assuming that the barriers we deal with consist
exclusively of singletons (an array with domain such a barrier is called 
a bad sequence because of the identification made in remark~\ref{wqo}) 
implies \PCA. This special case is known as the minimal bad 
sequence lemma and is a basic tool of wqo theory. Let us state it 
more explicitly. 

\begin{definition}
(\RCA) We say that a $Q$-sequence $\<q_n: n \in A\>$ with $A \in 
\Sq\N$ is bad if $\forall n,m \in A (n<m \imply q_n \npreceq 
q_m)$. If $<'$ is compatible with $\preceq$ and 
$\<q_n: n \in A\>$ and $\<q'_n: n \in A'\>$ are two $Q$-sequences  
we write $\<q_n\> \leq' \<q'_n\>$ if $A \subseteq A'$ and 
$\forall n \in A\; q_n \leq' q'_n$; we write $\<q_n\> <' 
\<q'_n\>$ if $\<q_n\> \leq' \<q'_n\>$ and $\exists n \in A\; q_n 
<' q'_n$. $\<q_n\>$ is minimal bad with respect to $<'$ if it is 
bad and there is no $\<q'_n\> <' \<q_n\>$ which is bad. 
\end{definition}

\begin{definition}
The minimal bad sequence lemma is the following statement: {\em 
let $(Q, \preceq)$ be a quasi-ordering and $<'$ a well-founded 
relation on $Q$ which is compatible with $\preceq$: if there 
exists a bad $Q$-sequence then there exists a $Q$-sequence which 
is minimal bad with respect to $<'$.}
\end{definition} 

We start our proof that the minimal bad sequence lemma implies 
\PCA\ by deducing \ACA.

\begin{lemma}\label{mbsACA}
(\RCA) The minimal bad sequence lemma implies \ACA.
\end{lemma}
\begin{pf}
It is well-known that \ACA\ is equivalent to \S01 comprehension
over \RCA. Let $\psi(n,m)$ be a \D00 formula. Reasoning in \RCA\ 
we want to show 
that there exists a set $Z$ such that $\forall n (n \in Z \sse 
\exists m \psi (n,m))$. Let $Q$ be the set of all finite 
sequences of natural numbers and define
$$s \preceq t \sse \lh(s) = \lh(t) \land \forall n < \lh(s) 
(\psi(n,t(n)) \imply \psi(n,s(n)))$$
Let $s <' t$ if and only if $s \preceq t$ and $t \npreceq s$.

$<'$ is well-founded on $Q$ and there exists a bad $Q$-sequence 
with respect to $\preceq$: we are therefore in the hypothesis of 
the minimal bad sequence lemma. Let $\<s_k: k \in A\>$ be a minimal bad
$Q$-sequence. Let $n_k = \lh(s_k)$ for every $k$: since $\<s_k\>$ is 
bad it is immediate that for every $n$ there exists only finitely many
$k \in A$ such that $n_k = n$. We claim that
$$\forall n \left( \exists m \psi(n,m) \sse \forall k \in A (n_k > n 
\imply \psi(n, s_k(n)))\right)$$
One direction follows from the above observation. For the other 
direction suppose $\psi(n,m)$ holds but $i$ is such that $n_i > n$ 
and $\psi(n, s_i(n))$ does not hold. Let $t_i$ be such that $\lh(t_i) 
= n_i$, $t_i(n) = m$ and $t_i(j) = s_i(j)$ if $j \neq n$. Then $t_i
<' s_i$. Let $t_k = s_k$ if $k \neq i$ and $A' = \set{k \in A}{k = i 
\lor n_k \neq n_i}$. $\<t_k : k \in A'\>$ is bad and $\<t_k : k 
\in A'\> <' \<s_k : k \in A\>$, contradicting the minimality of 
$\<s_k\>$.

The set $Z$ can now be defined by \D01 comprehension, which is 
available in \RCA.
\end{pf}

The following definition is helpful in the proof of the main 
result of this section. 
   
\begin{definition}
(\RCA) If $T$ is a tree and $f:\N \to \N$ is a path through $T$ 
(i.e.\ for every $n$ we have $f[n] \in T$) we say that $f$ is the 
leftmost path through $T$ if for every path $g \neq f$ through 
$T$ there exists $n$ such that $f[n] = g[n]$ and $f(n) < g(n)$.
\end{definition}

\begin{theorem}
(\RCA) The following are equivalent:\begin{enumerate}
\item \PCA;
\item The locally minimal bad array lemma;
\item The minimal bad sequence lemma;
\item If $T$ is a non well-founded tree then there exists the leftmost 
path through $T$.\end{enumerate}
\end{theorem}
\begin{pf}
(1) implies (2) is theorem~\ref{lmba} and (3) is a special case of 
(2): hence it suffices to prove that (3) implies (4) and that (4) 
implies (1).

Assume the minimal bad sequence lemma: by lemma~\ref{mbsACA} we 
can reason in \ACA. Let $T$ be a non well-founded tree. Let $Q = T$ 
and define $\preceq$ to be the Kleene-Brouwer order $\KB$. Recall 
that $s \KB t$ if and only if either $t \initeq s$ or for some $i 
< \lh(s), \lh(t)$ we have $s[i] = t[i]$ and $s(i) < t(i)$. $\KB$ 
is a linear order and, since $T$ is not well-founded, is not a 
well-ordering on $T$ (this can easily be proved in \RCA) and therefore 
is not a wqo. Define $<'$ by $s <' t$ if and only if $s <_{KB} t \land 
\lh(s) \leq \lh(t)$. $<'$ is well-founded on $T$ because for every $t 
\in T$ the set $\set{s}{s <' t}$ is a well-founded tree and hence 
$\KB$, and therefore also $<'$, is well-founded on it (here reasoning 
in \RCA\ would not suffice). 

Thus we are in the hypothesis of the minimal bad sequence lemma and
there exists a minimal bad $T$-sequence $\<s_n\>$ with respect to 
$<'$. Notice that $\<s_n\>$ is actually a descending chain with 
respect to $\KB$. For notational simplicity we  assume that $s_n$ 
is defined for every $n$. We claim that for every $n$ $s_n \init 
s_{n+1}$. Suppose this is not the case and $s_m \not\init s_{m+1}$: 
since $s_{m+1} <_{KB} s_m$ there exists $i < \lh(s_m), \lh(s_{m+1})$
such that $s_m[i] = s_{m+1}[i]$ and $s_m(i) > s_{m+1}(i)$. Let $t_m 
= s_{m+1}[i+1]$, so that $t_m <' s_m$ and $s_{m+1} \KB t_m$. Moreover 
let $t_n = s_n$ if $n \neq m$ and $A = \{m\} \cup \set{n}{n > m+1}$. 
Then $\<t_n: n \in A\>$ is bad and $\<t_n: n \in A\> <' \<s_n: n \in 
\N\>$, a contradiction which establishes the claim. 

The claim implies that there exists $f: \N \to \N$ such that for 
every $n$ $s_n \init f$. It is easy to check that $f$ is the leftmost 
path through $T$.

To show that (4) implies (1) recall that \PCA\ is equivalent 
(over \RCA) to the following statement: if $\<T_n\>$ is a 
sequence of trees there exists a set $Z$ such that $\forall n (n 
\in Z \sse T_n \text{ is well-founded})$. Let $\<T_n\>$ be given. 
For every $n$ let $T'_n = \set{\<0\> \conc s}{s \in T_n} \cup 
\set{s}{\forall i < \lh(s)\; s(i) =1}$. Clearly $T'_n$ is 
not well-founded. Let $T$ be the tree obtained by interleaving 
the $T'_n$'s, i.e.
$$T = \set{t}{\forall (i,n) < \lh(t)\; \<t((0,n)), \dots 
t((i,n))\> \in T'_n}$$
$T$ is not well-founded and hence by (4) it has a leftmost path 
$f$. Let $Z = \set{n}{f((0,n))=1}$. It is easy to check that $Z$ 
has the desired property. \end{pf} 

The binary relation used in proving \GHT\ in lemma~\ref{GHT} has
the property that for every $q \in Q$ the set $\set{q' \in Q}{q'
<' q}$ is finite. We do not know whether the minimal bad sequence 
lemma restricted to $<'$ having this additional property implies \PCA.

We would like to mention also that a more general version of the 
locally minimal bad array lemma is provable in \PCA. This version 
asserts that the lemma holds also if $Q$ is an (uncountable) 
analytic set, $\preceq$ is coanalytic and $<'$ is analytic: the 
proof consists of carrying out part of the proof of 
theorem~\ref{lmba} inside a $\beta$-model.

\end{document}